\documentstyle[amssymb,thmsc,sw20lart]{article}

\input tcilatex
\QQQ{Language}{
American English
}

\begin{document}

\author{Alexander Dvorsky and Siddhartha Sahi\thinspace \thanks{%
\ E-mail: {\tt dvorsky@math.rutgers.edu, \ sahi@math.rutgers.edu}}\smallskip 
\\
Department of Mathematics, Rutgers University, New Brunswick\\
NJ 08903, USA}
\title{Explicit Hilbert spaces for certain unipotent representations II.}
\date{}
\maketitle

\section{Introduction}

To each real semisimple Jordan algebra, the Tits-Koecher-Kantor theory
associates a distinguished parabolic subgroup $P=LN$ of a semisimple Lie
group $G$. The groups $P$ which arise in this manner are precisely those for
which $N$ is abelian, and $P$ is conjugate to its opposite $\overline{P}.$

Each non-open $L$-orbit ${\cal O}$ on $N^{*}$ admits an $L$-equivariant
measure $d\mu $ which is unique up to scalar multiple. By Mackey theory, we
obtain a natural irreducible unitary representation $\pi _{{\cal O}}$ of $P$%
, acting on the Hilbert space 
\[
{\cal H}_{{\cal O}}=L^2({\cal O},d\mu ). 
\]

In this context, we wish to consider two problems:

\begin{enumerate}
\item  Extend $\pi _{{\cal O}}$ to a unitary representation of $G.$

\item  Decompose the tensor products $\pi _{{\cal O}}\otimes \pi _{{\cal O}%
^{\prime }}\otimes \pi _{{\cal O}^{\prime \prime }}\otimes \cdots $
\end{enumerate}

If the Jordan algebra is Euclidean (i.e. formally real) then $G/P$ is the
Shilov boundary of a symmetric tube domain. In this case, the first problem
was solved in \cite{sahi-expl}, \cite{shilov}, where it was shown that $\pi
_{{\cal O}}$ extends to a unitary representation of a suitable covering
group of $G$. The second problem was solved in \cite{tens}, where we
established a correspondence between the unitary representations of $G$
occurring in the tensor product, and those of a ``dual'' group $G^{\prime }$
acting on a certain reductive homogeneous space. This correspondence agrees
with the $\theta $-correspondence in various classical cases, and also gives
a duality between $E_7$ and real forms of the Cayley projective plane.

In this paper we start to consider these two problems for {\em non-Euclidean}
Jordan algebras. The algebraic groundwork has already been accomplished in 
\cite{sahi-dp}, however the analytical considerations are much more subtle,
and here we only treat the case of the representation $\pi _1=\pi _{{\cal O}%
_1}$ corresponding to the {\em minimal} $L$-orbit ${\cal O}_1$.

It turns out that in order for the first problem to have a positive
solution, one has to exclude certain Jordan algebras of rank $2.$ This is
related to the Howe-Vogan result on the non-existence of minimal
representations for certain orthogonal groups.

To each of the remaining Jordan algebras we attach a restricted root system $%
\Sigma $ of rank $n$, where $n$ is the rank of the Jordan algebra. The root
multiplicities, $d$ and $e$, of $\Sigma $ play a decisive role in our
considerations. For the reader's convenience, we include a list of the
corresponding groups $G$ and the multiplicities in the appendix.

For these groups, we show that $\pi _1$ extends to a spherical unitary
representation of $G$, and that the spherical vector is closely related to
the {\em one}{\it \/} variable Bessel $K$-function $K_\tau (z)$, where

\[
\tau =\frac{d-e-1}2. 
\]

The function $K_\tau (z)$ can be characterized, up to a multiple, as the
unique solution of the modified Bessel equation 
\[
\psi ^{\prime \prime }+z^{-1}\psi ^{\prime }-(1+\frac{\tau ^2}{z^2})\psi =0 
\]
that decays (exponentially) as $z\rightarrow \infty $; and, to us, one of
the most delightful aspects of the present consideration is the unexpected
and uniform manner in which this classical differential equation emerges
from the structure theory of $G$.

More precisely, we establish the following result:

We identify $N$ with its Lie algebra ${\frak n}=\limfunc{Lie}(N)$ via the
exponential map. We also fix an invariant bilinear form on $\left\langle
\cdot ,\cdot \right\rangle $ on ${\frak g}$, which is a certain multiple of
the Killing form, normalized as in Definition \ref{=form} below. We use this
form to identify $N^{*}$ with $\overline{{\frak n}}=\limfunc{Lie}(\overline{N%
})$. For $y$ in $\overline{{\frak n}}$, $\left\langle -\theta
y,y\right\rangle $ is positive, and we define 
\[
|y|=\sqrt{\left\langle -\theta y,y\right\rangle }. 
\]

\begin{theorem}
$\pi _1$ extends to a unitary representation of $G$ with spherical vector $%
|y|^{-\tau }K_\tau (|y|)$.\label{=Main-theorem}
\end{theorem}

Since $\pi _1$ is spherical, its Langlands parameter is its infinitesimal
character, and this can be determined via the (degenerate) principal series
imbedding described in section 2 below. It is then straightforward to verify
that $\pi _1$ is the minimal representation of $G$, with annihilator equal
to the Joseph ideal. (For $G=GL(n),$ the minimal representation is not
unique.)

Thus our construction should be compared to other realizations of the
minimal representations in \cite{brylinski}, \cite{torasso}, \cite{huang}
etc. Although our construction is for a more restrictive class of groups, it
does offer two advantages over the other constructions. The first advantage
is that our construction works for a larger class of representations, and
the second advantage is that it is well-suited for tensor product
computations.

Both of these features will be explored in detail in a subsequent paper. In
the present paper, we consider $k$-fold tensor powers of $\pi _1$, where $k$
is strictly smaller than $n$ (rank of $\Sigma $), and show that the
decomposition can be understood in terms of certain reductive homogeneous
spaces 
\[
G_k/H_{k\text{ }},1<k<n. 
\]
These spaces are defined in section 3, and are listed in the appendix.

We consider also the corresponding Plancherel decomposition: 
\[
L^2(G_k/H_k)=\int_{\widehat{G}_k}^{\oplus }m(\kappa )\kappa \,d\mu (\kappa
)\,\text{,} 
\]
where $d\mu $ is the Plancherel measure, and $m(\kappa )$ is the
multiplicity function. Then we have

\begin{theorem}
For $1<k<n$, there is a correspondence $\theta _k$ between $\widehat{G}_k$
and $\widehat{G},$ such that 
\[
\pi _1^{\otimes k}=\int_{\widehat{G}_k}^{\oplus }m(\kappa )\theta _k\left(
\kappa \right) d\mu (\kappa ).
\]
\label{=Main-tensoring}
\end{theorem}

\section{Preliminaries}

The results of this section are all well-known. Details and proofs may be
found in \cite{sahi-expl}, \cite{kostant-sahi} and in the references therein
(in particular, \cite{braun} and \cite{loos}).

\subsection{Root multiplicities}

Let $G$ be a real simple Lie group and let $K$ be a maximal compact subgroup
corresponding to a Cartan involution $\theta .$ We shall denote the Lie
algebras of $G$, $K$ etc by ${\frak g}$, ${\frak k}$ etc. Their
complexifications will be denoted by lowercase fraktur letters with
subscript $_{{\Bbb C}}$. Fix $\theta $, and let ${\frak g}={\frak k}+{\frak p%
}$ be the associated Cartan decomposition.

The parabolic subgroups $P=LN$ obtained by the Tits-Kantor-Koecher
construction are those such that $N$ is abelian, and $P$ is $G$-conjugate to
its opposite parabolic 
\[
\overline{P}=\theta (P)=L\overline{N}. 
\]
In this case $N$ has a natural structure of a real Jordan algebra, which is
unique up to a choice of the identity element.

In (Lie-)algebraic terms, this means that $P$ is a maximal parabolic
subgroup corresponding to a simple\ (restricted) root $\alpha $ which has
coefficient $1$ in the highest root, and which is mapped to $-\alpha $ under
the long element of the Weyl group.

In this situation, $M:=K\cap L$ is a symmetric subgroup of $K$ (this is {\em %
equivalent} to the abelianness of $N),$ and we fix a maximal toral
subalgebra ${\frak t}$ in the orthogonal complement of ${\frak m}$ in $%
{\frak k}$.

The roots of ${\frak t}_{{\Bbb C}}$ in ${\frak g}_{{\Bbb C}}$ form a
restricted root system of type $C_n$, where $n=\dim _{{\Bbb R}}{\frak t}$ is
the (real) rank of $N$ as a Jordan algebra (this result is essentially due
to C. Moore). We fix a basis $\{\gamma _1,\gamma _2,\ldots ,\gamma _n\}$ of $%
{\frak t}^{*}$ such that 
\[
\Sigma ({\frak t}_{{\Bbb C}},{\frak g}_{{\Bbb C}})=\{\pm (\gamma _i\pm
\gamma _j)/2,\pm \gamma _j\}\text{.} 
\]

The restricted root system $\Sigma =\Sigma ({\frak t}_{{\Bbb C}},{\frak k}_{%
{\Bbb C}})$ is of type $A_{n-1},C_n$ or $D_n$, and the first of these cases
arises precisely when $N$ is a Euclidean Jordan algebra. This case was
studied in \cite{sahi-expl}, therefore we restrict our attention to the last
two cases.

The root multiplicities in $\Sigma $ play a key role in our considerations.
If $\Sigma $ is $C_n$, there are two multiplicities, corresponding to the
short and long roots, which we denote by $d$ and $e,$ respectively. If $%
\Sigma $ is $D_n$, and $n\neq 2,$ then there is a single multiplicity, which
we denote by $d$, so that $D_n$ may be regarded as a special case of $C_n$,
with $e=0.$

The root system $D_2$ is reducible (being isomorphic to $A_1\times A_1$) and 
{\em a priori} there are two root multiplicities. In what follows, we
explicitly{\em \ }exclude the case when these multiplicities are different.%
{\em \ }This means that we exclude from consideration the groups 
\[
G=O(p,q),N={\Bbb R}^{p-1,q-1}(p\neq q); 
\]
indeed, our main results are false for these groups.{\em \ }When the two
multiplicities {\em coincide}, we once again denote the common multiplicity
by $d$.

The multiplicity of the short roots $\pm (\gamma _i\pm \gamma _j)/2$ in $%
\sum ({\frak t}_{{\Bbb C}},{\frak g}_{{\Bbb C}})$ is equal to $2d$, and the
multiplicity of the long roots $\pm \gamma _i$ is $e+1$.

In the appendix we include a table listing the groups under consideration,
as well as the values of $d$ and $e$ for each of these groups.

\subsection{Cayley transform}

We briefly review the notion of the Cayley transform. Let $C$ be the
following element (of order $8$) in $SL_2\left( {\Bbb C}\right) $ 
\[
C=\frac 1{\sqrt{2}}\left( 
\begin{array}{cc}
1 & i \\ 
i & 1
\end{array}
\right) . 
\]

The Cayley transform of ${\frak sl}_2({\Bbb C})$ is the automorphism (of
order $4$) given by 
\[
c=\text{Ad }C. 
\]

It transforms the ``usual'' basis of ${\frak sl}_2({\Bbb C})$ 
\[
x=\left( 
\begin{array}{cc}
0 & 1 \\ 
0 & 0
\end{array}
\right) ,\text{ }y=\ \left( 
\begin{array}{cc}
0 & 0 \\ 
1 & 0
\end{array}
\right) ,\text{ }h=\left( 
\begin{array}{cc}
1 & 0 \\ 
0 & -1
\end{array}
\right) , 
\]
to the basis

\[
X=\frac 12\left( 
\begin{array}{cc}
-i & 1 \\ 
1 & i
\end{array}
\right) ,\text{ }Y=\frac 12\left( 
\begin{array}{cc}
i & 1 \\ 
1 & -i
\end{array}
\right) ,\text{ }H=i\left( 
\begin{array}{cc}
0 & -1 \\ 
1 & 0
\end{array}
\right) , 
\]
where $X=c(x)=C^{-1}xC,$ etc. In turn, $c$ can be expressed as 
\[
c=\exp \text{ ad }\frac{\pi i}4(x+y)=\exp \text{ ad }\frac{\pi i}4(X+Y). 
\]

The key property of the Cayley transform is that it takes the compact torus
(spanned by $iH)$ to the split torus spanned by $h$ (cf. \cite{koranyi-wolf}%
).

We turn now to the Lie algebra ${\frak g}_{{\Bbb C}}$. By the
Cartan-Helgason theorem the root spaces ${\frak p}_{\gamma _j}$ are
one-dimensional, and so by the Jacobson-Morozov theorem we get holomorphic
homomorphisms 
\[
\Phi _j:{\frak sl}_2({\Bbb C})\longrightarrow {\frak g}_{{\Bbb C}}\text{, }%
j=1,...,n 
\]
such that $X_j=\Phi _j(X)$ spans ${\frak p}_{\gamma _j}$.

We fix such maps $\Phi _j$, and denote the images of $x,X,y,Y,h,H$ by $%
x_j,X_j,$ etc. Since the roots $\gamma _j$ are strongly orthogonal, the
triples $\{X_j,Y_j,H_j\}$ commute with each other, and the Cayley transform
of ${\frak g}$ is defined to be the automorphism 
\[
c=\exp \text{ ad }\frac{\pi i}4\left( \sum X_j+\sum Y_j\right) =\prod \exp 
\text{ ad }\frac{\pi i}4(X_j+Y_j). 
\]

Thus we obtain an ${\Bbb R}$-split toral subalgebra ${\frak a}$ defined by 
\[
{\frak a}=c^{-1}(i{\frak t})={\Bbb R}h_1\oplus \cdots \oplus {\Bbb R}h_n. 
\]
The roots of ${\frak a}_{{\Bbb C}}$ in ${\frak g}_{{\Bbb C}}$ are 
\[
\Sigma ({\frak a}_{{\Bbb C}},{\frak g}_{{\Bbb C}})=\left\{ \pm \varepsilon
_i\pm \varepsilon _j,\pm 2\varepsilon _j\right\} \text{ where }\varepsilon
_i=\frac 12\gamma _i\ \circ c. 
\]
The short roots have multiplicity $2d$ and the long roots have multiplicity $%
e+1.$

In fact ${\frak a}\subset {\frak l},$ and we have 
\[
\Sigma ({\frak a},{\frak l})=\left\{ \pm (\varepsilon _i-\varepsilon
_j)\right\} ,\text{ }\Sigma ({\frak a},{\frak n})=\left\{ \varepsilon
_i+\varepsilon _j,2\varepsilon _j\right\} ,\text{ }\Sigma ({\frak a},%
\overline{{\frak n}})=\left\{ -\varepsilon _i-\varepsilon _j,-2\varepsilon
_j\right\} 
\]

\begin{definition}
The invariant form $\left\langle .,.\right\rangle $ on ${\frak g}$ is
normalized by requiring 
\[
\left\langle x_1,y_1\right\rangle =1.
\]
\label{=form}
\end{definition}

For $y\in \overline{{\frak n}}$, we set $\left| y\right| \stackrel{\text{def}%
}{=}\sqrt{-\left\langle y,\theta y\right\rangle }$ , as in Introduction.

\subsection{ Orbits and measures}

We now describe the orbits of $L$ in $\overline{{\frak n}}$ $\simeq N^{*}.$
For $k=1,...,n-1,$ define 
\[
{\cal O}_k=L\cdot (y_1+y_2+\ldots +y_k). 
\]
Then these, together with the trivial orbit ${\cal O}_0$, comprise the
totality of the singular (i.e., non-open) $L$-orbits in $\overline{{\frak n}}
$.

We define $\nu \in {\frak a}^{*}$ as 
\[
\nu =\varepsilon _1+\varepsilon _2+\ldots +\varepsilon _n.\text{ } 
\]
Then $\nu $ extends to a character of ${\frak l}$ , and we will write $e^\nu 
$ for the corresponding (spherical) character of $L.$

\begin{lemma}
The orbit ${\cal O}_1$ carries a natural $L-$equivariant measure $d\mu _1$,
which transforms by the character $e^{2d\nu }$, that is 
\[
\int_{{\cal O}_1}g(l\cdot y)d\mu _1(y)=e^{2d\nu }(l)\int_{{\cal O}%
_1}g(y)d\mu _1(y)\text{.}
\]
\end{lemma}

\TeXButton{Proof}{\proof} Let $S_1$ be the stabilizer of $y_1$ in $L$. It
suffices to show that the modular function of$\,$ $S_1$ is the restriction,
from $L$ to $S_1$, of the character $e^{2d\nu }$. Passing to the Lie algebra 
${\frak s}_1$, we need to show that 
\[
\limfunc{tr}\limfunc{ad}\nolimits_{{\frak s}_1}=2d\nu |_{{\frak s}_1}. 
\]

To see this, we remark that ${\frak s}_1$ has codimension $1$ inside a
maximal parabolic subalgebra ${\frak q}$ of ${\frak l}$, corresponding to
the stabilizer of the line through $y_1.$ The space of characters of ${\frak %
q}$ is two-dimensional, and it follows that the space of characters of $%
{\frak s}_1$ is one-dimensional. Hence any character of ${\frak s}_1$ is
determined by its restriction to ${\frak a}\cap {\frak s}_1=\limfunc{Ker}%
\varepsilon _1$. The restriction of $\nu $ to ${\frak s}_1$ is nontrivial,
hence 
\[
\limfunc{tr}\limfunc{ad}\nolimits_{{\frak s}_1}=k\nu 
\]
for some constant $k$.

Obviously, $\limfunc{tr}\limfunc{ad}_{{\frak l}}=0$, and the only root
spaces missing from ${\frak s}_1$ are the root spaces ${\frak l}%
_{\varepsilon _1-\varepsilon _j}$, $j\geq 2$ (each of these root spaces has
dimension $2d$). Hence, for $a\in {\frak a}$%
\[
\limfunc{tr}\limfunc{ad}\nolimits_{{\frak s}_1}(a)=-2d\sum_{j=2}^n(%
\varepsilon _1-\varepsilon _j)(a)\text{,} 
\]
and restricting this to $\limfunc{Ker}\varepsilon _1$, we obtain $2d\nu \,|_{%
{\frak a}\cap {\frak s}_1}$.\TeXButton{End Proof}{\endproof}

\smallskip\ 

{\em Example.} Consider $G=O_{2n,2n}$ realized as the group of all $2n\times
2n$ real matrices preserving the split symmetric form $\left( 
\begin{array}{cc}
0 & I_{2n} \\ 
I_{2n} & 0
\end{array}
\right) $. Then $P=LN=GL_{2n}({\Bbb R})\rightthreetimes \limfunc{Skew}_{2n}(%
{\Bbb R}).$ More precisely, 
\[
L=\left\{ \left( 
\begin{array}{cc}
A & 0 \\ 
0 & A^{t^{-1}}
\end{array}
\right) :A\in GL_{2n}({\Bbb R})\right\} 
\]
and 
\[
N=\left\{ \left( 
\begin{array}{cc}
I_{2n} & 0 \\ 
B & I_{2n}
\end{array}
\right) :B+B^t=0\right\} . 
\]
Then 
\[
{\frak a}=\left\{ \limfunc{diag}(a_1,a_1,a_2,a_2,\ldots
,a_n,a_n,-a_1,-a_1,-a_2,-a_2,\ldots ,-a_n,-a_n),\,a_i\in {\Bbb R}\right\} 
\]
is the toral subalgebra of ${\frak g}$ (and ${\frak l}$) described in the
preceding subsection. We can take 
\[
y_1=\left( 
\begin{array}{cc}
0_{2n} & B_1 \\ 
0 & 0_{2n}
\end{array}
\right) \text{, where }B_1=\left( 
\begin{array}{ccc}
0 & -1 & 0 \\ 
1 & 0 & 0 \\ 
0 & 0 & 0
\end{array}
\right) . 
\]
The Lie algebra ${\frak s}_1$ of the stabilizer $S_1=\limfunc{Stab}_Ly_1$
can be written as 
\[
{\frak s}_1=\left\{ \left( 
\begin{array}{cc}
A & 0 \\ 
0 & -A^t
\end{array}
\right) :A=\left( 
\begin{array}{cc}
A_{11} & 0 \\ 
A_{21} & A_{22}
\end{array}
\right) ,A_{11}\in {\frak sl}_2,A_{22}\in {\frak gl}_{2n-2}\right\} 
\]
It is a codimension 1 subalgebra of the parabolic subalgebra ${\frak q}$ of $%
{\frak gl}_{2n}$, where ${\frak q}=({\frak gl}_2+{\frak gl}_{2n-2})+{\Bbb R}%
^{2,2n-2}.$\medskip\ 

{\em Remark.} In this example $\nu =\frac 12\limfunc{tr}$, $d=2$ and $%
e^{2d\nu }=(\det )^2$.

\section{Minimal representation of $G$}

If $\chi $ is a character of ${\frak l}$, we write $\pi _\chi $ for the
(unnormalized) induced representation $\limfunc{Ind}_{\overline{P}}^G(\chi )$%
. These representations were studied in \cite{sahi-dp} in the ``compact''
picture, by algebraic methods. Among the results established there was the
existence of a finite number of ``small'', unitarizable, spherical
subrepresentations, which occur for the following values of $\chi $ 
\[
\chi _j=e^{-jd\nu },\text{ }j=1,\ldots ,n-1. 
\]

In this paper we use analytical methods, and work primarily with the
``non-compact'' picture, which is the realization of $\pi _\chi $ on $%
C^\infty (N)$, via the Gelfand-Naimark decomposition

\[
G\approx N\overline{P}. 
\]
In fact, using the exponential map we can identify ${\frak n}$ and $N$, and
realize $\pi _\chi $ on $C^\infty ({\frak n}).$

We will show that the unitarizable subrepresentation of $\pi _{\chi _1}$
admits a natural realization on the Hilbert space $L^2({\cal O}_1,d\mu )$.
Since there is no obvious action of $G$ on this space, we have to proceed in
an indirect fashion. The key is an explicit realization of the spherical
vector $\sigma _{\chi _1}.$

\subsection{The Bessel function}

We let $d,e$ be the root multiplicities of $\Sigma ({\frak t},{\frak k})$ as
in previous section, and define

\[
\tau _G=\tau =(d-e-1)/2 
\]
as in the introduction.

Let $K_\tau $ be the $K$-Bessel function on $(0,\infty )$ satisfying

\begin{equation}
z^2K_\tau ^{\prime \prime }+zK_\tau ^{\prime }-(z^2+\tau ^2)K_\tau =0.
\label{=bessel-eq}
\end{equation}

Put $\phi _\tau (z)=\dfrac{K_\tau (\sqrt{z})}{\left( \sqrt{z}\right) ^\tau }$%
, then $\phi _\tau $ satisfies the differential equation

\begin{equation}
D\phi _\tau =0,\text{where }D\phi =4z\phi ^{\prime \prime }+4(\tau +1)\phi
^{\prime }-\phi \text{.}  \label{=change}
\end{equation}

We lift $\phi _\tau $ to an $M$-invariant function $g_\tau $ on ${\cal O}_1$%
, by defining

\begin{equation}
g_\tau (y)=\phi _\tau (-\left\langle y,\theta y\right\rangle )=\frac{K_\tau
\left( \left| y\right| \right) }{\left| y\right| ^\tau }.  \label{=g-tau}
\end{equation}

{\em Remark}. If $d=e$ (as is the case for $G=Sp_{2n}({\Bbb C})$ or $%
Sp_{n,n} $), then $\tau =-\frac 12$ and 
\[
g_\tau (y)=\left| y\right| ^{1/2}K_{-1/2}(\left| y\right| )=\left| y\right|
^{1/2}\frac{\exp (-\left| y\right| )}{\left| y\right| ^{1/2}}=e^{-\left|
y\right| }. 
\]
If $d=e+1$ (this is true for $GL_{2n}({\bf k})$, ${\bf k}={\Bbb R}$, ${\Bbb C%
}$ or ${\Bbb H}$), then 
\[
g_\tau (y)=K_0(\left| y\right| ). 
\]

\begin{proposition}
{\em (1)} $g_\tau $ is a (square-integrable) function in $L^2({\cal O}%
_1,d\mu _1)$.

{\em (2) }The measure $g_\tau d\mu _1$ defines a tempered distribution on $%
\overline{{\frak n}}$.\label{=sq-int}
\end{proposition}

\TeXButton{Proof}{\proof} (1) We define 
\[
{\cal O}^{\prime }\stackrel{\text{def}}{=}\{y^{\prime }\in {\cal O}_1:\left|
y^{\prime }\right| =1\}. 
\]
Then ${\cal O}^{\prime }$ is compact; the map 
\[
{\cal O}^{\prime }\times (0,\infty )\ni (y^{\prime },w)\longmapsto
wy^{\prime }\in {\cal O}_1 
\]
is a diffeomorphism, and the measure $d\mu _1$ can be decomposed as a product

\[
d\mu _1(wy^{\prime })=d\mu ^{\prime }(y^{\prime })d\mu ^{\prime \prime }(w) 
\]

We now determine the explicit form of $d\mu ^{\prime \prime }(w).$

Define $h=\sum_{i=1}^nh_i$, then $(\limfunc{ad}h)y=-2y$ for any $y\in 
\overline{{\frak n}}$. We take $y\in {\cal O}_1$, $z>0$, $a=\ln z$ and
calculate 
\[
d\mu _1(zy)=d\mu _1(\exp (-a\frac h2)\cdot y)=e^{-2d\nu (-\frac{ah}2)}d\mu
_1(y)=e^{dna}d\mu _1(y)=z^{dn}d\mu _1(y)\text{.} 
\]
Therefore, for $z>0$ 
\begin{equation}
d\mu _1(zy)=z^{dn}d\mu _1(y)  \label{=scaling}
\end{equation}
and it follows that $d\mu ^{\prime \prime }(zw)=z^{dn}d\mu ^{\prime \prime
}(w),$ and so, up to a scalar multiple, 
\[
d\mu ^{\prime \prime }(w)=w^{dn-1}dw, 
\]
where $dw$ is the Lebesgue measure.

We can now calculate 
\begin{eqnarray}
\int_{{\cal O}_1}\left| g_\tau (y)\right| ^2d\mu _1(y) &=&\int_0^\infty
\int_{{\cal O}^{\prime }}\frac{K_\tau (w)^2}{w^{2\tau }}d\mu ^{\prime
}(y^{\prime })w^{dn-1}dw  \nonumber \\
&=&c\int_0^\infty \frac{K_\tau (w)^2}{w^{2\tau }}w^{dn-1}dw,
\label{=sep-var}
\end{eqnarray}
where $c=\mu ^{\prime }({\cal O}^{\prime })$ is a positive constant. The
function $K_\tau (w)$ has a pole of order $\tau $ at $0$ (or, in case of $%
\tau =0,$ a logarithmic singularity at $0$), and it decays exponentially as $%
w\rightarrow \infty $ \cite[3.71.15]{watson}. Hence $w^{-2\tau }K_\tau (w)^2$
has a pole of order 
\[
4\tau =2(d-e-1)\leq 2d-2<dn-1 
\]
(recall that we require $n\geq 2)$. Thus the integrand in (\ref{=sep-var})
is non-singular and decays exponentially as $w\rightarrow \infty $.
Therefore, the integral (\ref{=sep-var}) converges and $g_\tau (y)\in L^2(%
{\cal O}_1,d\mu _1)$.

(2) From the calculation in (1), we see that $g_\tau (y)\in L_{\text{loc}}^1(%
{\cal O}_1,d\mu _1)$ and has exponential decay at $\infty $ (i.e., as $%
\left| y\right| \rightarrow \infty $). This implies the result. 
\TeXButton{End Proof}{\endproof}\smallskip

We can now define the Fourier transform of $g_\tau $, 
\[
\Phi =\widehat{g_\tau d\mu _1} 
\]
as a (tempered) distribution on ${\frak n.}$ The key result is the following

\begin{proposition}
$\Phi $ is a multiple of the spherical vector $\sigma _{\chi _1}$.\label
{=spher}
\end{proposition}

The proof of this proposition will be given over the next two subsections.%
\label{=analytic}

\subsection{Characterization of spherical vectors}

For $\phi :{\frak n}\rightarrow {\frak n},$ let $\xi (\phi )$ denote the
corresponding vector field: 
\[
\xi (\phi )f(x)=\left. \frac d{dt}f(x+t\phi (x))\right| _{t=0}\text{ for }f:%
{\frak n}\rightarrow {\Bbb C}. 
\]

Then we have the following formulas for the action of $\pi _\chi $ on $%
C^\infty ({\frak n})$:

\begin{itemize}
\item  for $x_0\in {\frak n}$, $\pi _\chi (x_0)=\xi (x_0)$,

\item  for $h_0\in {\frak l}$, $\pi _\chi (h_0)=\chi (h_0)-\xi \left(
[h_0,x]\right) $,

\item  for $y_0\in \overline{{\frak n}}$, $\pi _\chi (y_0)=\chi
[x,y_0]-\frac 12\xi \left( [h,x]\right) $, where $h=[x,y_0]$.
\end{itemize}

We need a Lie algebra characterization of $\sigma _\chi $:

\begin{lemma}
The space of $\pi _\chi ({\frak k})$-invariant distributions on ${\frak n}$
is 1--dimensional (and spanned by $\sigma _\chi $).
\end{lemma}

\TeXButton{Proof}{\proof} It is well known (and easy to prove) that the only
distributions on ${\Bbb R}^n,$ which are annihilated by $\frac \partial
{\partial x_i},$ $i=1,...,n$ are the constants. More generally, we can
replace ${\Bbb R}^n$ by a manifold, and $\left\{ \frac \partial {\partial
x_i}\right\} $ by any set of vector fields which span the tangent space at
each point of the manifold.

For $\chi =0$, the formulas above show that $\pi _0({\frak g})$ acts by
vector fields on $C^\infty ({\frak n}).$ Moreover, using the decomposition $%
G=K\overline{P},$ we see that $\pi _0({\frak k})$ is a spanning family of
vector fields. Thus the result follows in this case.

For general $\chi $, if $T$ is a $\pi _\chi ({\frak k})$-invariant
distribution, then $T/\sigma _\chi =T\sigma _{-\chi }$ is $\pi _0({\frak k})$%
-invariant, and hence a constant. \TeXButton{End Proof}{\endproof}

\begin{proposition}
Let $T$ be an $M$-invariant distribution on ${\frak n}$ such that 
\[
\pi _\chi (y+\theta y)T=0\,\text{ for {\em some} }y\text{ }\neq 0\text{ in }%
\overline{{\frak n}},
\]
then $T$ is a multiple of the spherical vector $\sigma _\chi $.\label{=infin}
\end{proposition}

\TeXButton{Proof}{\proof} The $M$-invariance of $T$ implies that 
\[
\pi _\chi ({\frak m})T=0 
\]
Since ${\frak m}$ is a maximal subalgebra of ${\frak k}$, ${\frak m}$ and $%
y+\theta y$ generate ${\frak k}$ as a Lie algebra. Thus 
\[
\pi _\chi ({\frak k})T=0, 
\]
and the result follows from the previous lemma. \TeXButton{End Proof}
{\endproof}

\subsection{The $K$-invariance of the Bessel function}

We now turn to the proof of Proposition \ref{=spher}. To simplify notation,
we will write $\pi $ instead of $\pi _{\chi _1}.$ Since $\Phi $ is clearly $%
M $-invariant, by Proposition \ref{=infin} it suffices to show 
\[
\pi (y_1+\theta y_1)\Phi =0 
\]
for $y_1\in \overline{{\frak n}}$. We will prove this through a sequence of
lemmas.

It is convenient to introduce the following notation: if $g_1$ and $g_2$ are
functions on ${\cal O}_1$, we define 
\[
(g_1,g_2)=\int_{{\cal O}_1}g_1(y)\,g_2(y)d\mu _1(y) 
\]
provided the integral converges.

If $g$ is a function on ${\cal O}_1$ and $h\in {\frak l}$, then the action
of $h$ on $g$ is given by 
\[
h\cdot g(y)\stackrel{\text{def}}{=}\left. \frac d{dt}g(e^{th}\cdot y)\right|
_{t=0}. 
\]

In the computation below, we shall work with the expressions of the type 
\[
\left. \left( \frac d{dt}\int_{{\cal O}_1}g(e^{th}\cdot y)d\mu (y)\right)
\right| _{t=0}. 
\]
To justify differentiation under the integral sign, we need to impose the
standard conditions on $g$ (e.g. \cite[p.170]{kestleman}), as follows.

Define a class of functions ${\cal I}\subset C^\infty ({\cal O}_1)$, given
by the following conditions: a smooth function $g\ $belongs to ${\cal I}$ if

\begin{itemize}
\item  $g\in L^1({\cal O}_1,d\mu _1)$ and

\item  for any $h\in {\frak l}$ we can find $c>0$ and $G(y)\in L^1({\cal O}%
_1,d\mu _1)$, such that 
\[
\left| \left. \frac d{dt}g(e^{th}\cdot y)\right| _{t=t_0}\right| \leq G(y)
\]
for all $y\in {\cal O}_1$ and $\left| t_0\right| <c$.
\end{itemize}

\begin{lemma}
Suppose $g_1,g_2$ are smooth functions on ${\cal O}_1$, such that $g_1g_2\in 
{\cal I}$. Then 
\begin{equation}
(h\cdot g_1,g_2)+(g_1,h\cdot g_2)=2d\nu (h)(g_1,g_2).  \label{=leib}
\end{equation}
\label{=leibnitz}
\end{lemma}

\TeXButton{Proof}{\proof} Using the $L$-equivariance of $d\mu _1$, we obtain 
\[
\int_{{\cal O}_1}g_1(e^{th}y)g_2(e^{th}y)\,d\mu _1=e^{2td\nu (h)}\int_{{\cal %
O}_1}g_1g_2\,d\mu _1. 
\]
Under the assumptions of the lemma, we can differentiate this identity in $t$%
, to get 
\[
\int_{{\cal O}_1}h\cdot (g_1g_2)\,d\mu _1=2d\nu (h)\int_{{\cal O}%
_1}g_1g_2\,d\mu _1\text{.} 
\]
By the Leibnitz rule, the result follows.\TeXButton{End Proof}{\endproof}

More generally, if $g_1,g_2$ are functions on ${\frak n}\times {\cal O}_1$,
then $(g_1,g_2)$ is a function on ${\frak n}$. In this notation, for $g$ in $%
L^1({\cal O}_1,d\mu _1)$, the Fourier transform of $gd\mu _1$ is given by
the formula 
\[
\widehat{gd\mu _1}=(e^{-i\left\langle x,y\right\rangle },g). 
\]

\begin{lemma}
Let $g\in L^1({\cal O}_1,d\mu _1)$ be a smooth function on ${\cal O}_1$,
such that 
\[
e^{-i\left\langle x,y\right\rangle }g\in {\cal I}.
\]
Suppose $f=(e^{-i\left\langle x,y\right\rangle },g)$, then 
\[
\pi (y_1)f=-\frac 12(e^{-i\left\langle x,y\right\rangle },h\cdot g(y))\text{%
, where }h=[x,y_1]\text{.\label{=pi-y0}}
\]
\end{lemma}

\TeXButton{Proof}{\proof} By the formula for the action of $\pi (y_1),$ we
get 
\begin{eqnarray*}
-2\left( \pi (y_1)f+d\nu (h)f\right) &=&\xi \left( [h,x]\right) \cdot
(e^{-i\left\langle x,y\right\rangle },g) \\
&=&\frac d{dt}\left. \left( e^{-i\left\langle x+t[h,x],y\right\rangle
},g\right) \right| _{t=0} \\
&=&\frac d{dt}\left. \left( e^{-i\left\langle x,y-t[h,y]\right\rangle
},g\right) \right| _{t=0} \\
&=&-\left( h\cdot e^{-i\left\langle x,y\right\rangle },g\right) \\
&=&(e^{-i\left\langle x,y\right\rangle },h\cdot g)-2d\nu
(h)(e^{-i\left\langle x,y\right\rangle },g).
\end{eqnarray*}
where we have used the previous lemma, and the relation 
\[
\left\langle x+t[h,x],y\right\rangle =\left\langle x,y\right\rangle
+t\left\langle [h,x],y\right\rangle =\left\langle x,y\right\rangle
-t\left\langle x,[h,y]\right\rangle =\left\langle x,y-t[h,y]\right\rangle . 
\]
The result follows.\TeXButton{End Proof}{\endproof}

The pairing $-\left\langle \cdot ,\theta \cdot \right\rangle $ gives a
positive definite $M$-invariant inner product on $\overline{{\frak n}}$, and
we now obtain the following

\begin{lemma}
Suppose that $g(y)=\phi \left( -\left\langle y,\theta y\right\rangle \right) 
$ for some smooth $\phi $ on $\left( 0,\infty \right) $, and $%
e^{-i\left\langle x,y\right\rangle }g\in {\cal I}$. Put $f(x)=(e^{-i\langle
x,y\rangle },g)$, as before. Then 
\[
\pi (y_1)f=\left( e^{-i\left\langle x,y\right\rangle },\left\langle
x,[[\theta y,y_1],y]\right\rangle \phi ^{\prime }\left( -\left\langle
y,\theta y\right\rangle \right) \right) .
\]
\end{lemma}

\TeXButton{Proof}{\proof} Writing $h=[x,y_1]$ as in the previous lemma, we
get 
\begin{eqnarray*}
h\cdot g(y) &=&\frac d{dt}\left. \phi \left( -\left\langle y+t[h,y],\theta
(y+t[h,y])\right\rangle \right) \right| _{t=0} \\
&=&\frac d{dt}\left. \phi \left( -\left\langle y,\theta y\right\rangle
-2t\left\langle \theta y,[h,y]\right\rangle +O(t^2)\right) \right| _{t=0} \\
&=&-2\left\langle \theta y,[h,y]\right\rangle \phi ^{\prime }\left(
-\left\langle y,\theta y\right\rangle \right) \text{.}
\end{eqnarray*}
Since 
\[
\left\langle \theta y,[h,y]\right\rangle =\left\langle \theta
y,[[x,y_1],y]\right\rangle =\left\langle x,[[\theta y,y_1],y]\right\rangle , 
\]
the result follows.\TeXButton{End Proof}{\endproof}

The key lemma is the following computation

\begin{lemma}
Let $\phi $ and $f$ be as in the previous lemma, and suppose for $x\in 
{\frak n}$ 
\begin{equation}
e^{-i\left\langle x,y\right\rangle }\phi \left( -\left\langle y,\theta
y\right\rangle \right) \in {\cal I}\text{, \ \ }e^{-i\left\langle
x,y\right\rangle }\phi ^{\prime }\left( -\left\langle y,\theta
y\right\rangle \right) \in {\cal I}\text{.}  \label{=diff-cond}
\end{equation}
Then we have 
\begin{equation}
\pi (y_1+\theta y_1)f(x)=\left( e^{-i\left\langle x,y\right\rangle
},i\left\langle \theta y_1,y\right\rangle \,(D\phi )\left( -\left\langle
y,\theta y\right\rangle \right) \right) ,  \label{=crown}
\end{equation}
where the differential operator $D$ is given by the formula {\em (\ref
{=change})}, i.e. 
\begin{equation}
(D\phi )\left( -\left\langle y,\theta y\right\rangle \right) =4\left(
-\left\langle y,\theta y\right\rangle \right) \phi ^{\prime \prime
}+2(d+1-e)\phi ^{\prime }-\phi \text{.}  \label{=de1}
\end{equation}
\label{=main-mess}
\end{lemma}

\TeXButton{Proof}{\proof} Choose a basis $l_j$ of ${\frak l}$ and define
functions $c_j(y)$ by the formula $[\theta y,y_1]=\sum_jc_j(y)l_j$. Then by
the previous lemma 
\begin{eqnarray*}
\pi (y_1)f &=&\sum_j(e^{-i\left\langle x,y\right\rangle },\left\langle
x,[l_j,y]\right\rangle c_j\phi ^{\prime })=i\sum_j\frac d{dt}\left. \left(
e^{-i\left\langle x,y+t[l_j,y]\right\rangle },c_j\phi ^{\prime }\right)
\right| _{t=0} \\
&=&i\sum_j(l_j\cdot e^{-i\left\langle x,y\right\rangle },c_j\phi ^{\prime })%
\text{.}
\end{eqnarray*}
Differentiation in this calculation is justified, because $e^{-i\left\langle
x,y\right\rangle }\phi ^{\prime }\left( -\left\langle y,\theta
y\right\rangle \right) \in {\cal I}$.

Applying (\ref{=leib}) to the last expression, we can write 
\begin{equation}
\pi (y_1)f=-i\sum_j\left( e^{-i\left\langle x,y\right\rangle },-2d\nu
(l_j)c_j\phi ^{\prime }+c_jl_j\cdot \phi ^{\prime }+\phi ^{\prime }l_j\cdot
c_j\right) \text{.}  \label{=pi-y1}
\end{equation}
We now calculate each term in this expression.

\begin{itemize}
\item  First we have 
\[
\sum_j\nu (l_j)c_j\phi ^{\prime }=\nu ([\theta y,y_1])\phi ^{\prime }.
\]
Since $\nu $ is a real character of ${\frak l}$, it vanishes on ${\frak l}%
\cap {\frak k}$ and we have $\nu ([\theta y,y_1])=\nu ([\theta y_1,y])$.
Recall that ${\frak n}$ and $\overline{{\frak n}}$ are irreducible ${\frak l}
$-modules. Therefore, $\nu ([\theta y_1,y])=k\left\langle \theta
y_1,y\right\rangle $ for some constant $k\neq 0,$ independent of $y$.
Setting $y=y_1$, we get $\left\langle \theta y_1,y_1\right\rangle
=\left\langle -x_1,y_1\right\rangle =-1$. Hence $k=-\nu ([\theta y_1,y_1])=1$%
, and therefore 
\begin{equation}
-\sum_j2d\nu (l_j)c_j\phi ^{\prime }=-2d\left\langle \theta
y_1,y\right\rangle \phi ^{\prime }\text{.}  \label{=term1}
\end{equation}

\item  Next we compute 
\begin{eqnarray*}
\sum_jc_jl_j\cdot \phi ^{\prime }=[\theta y,y_1]\cdot \phi ^{\prime } \\
=\frac d{dt}\left. \phi ^{\prime }\left( -\left\langle y+t[[\theta
y,y_1],y],\theta (y+t[[\theta y,y_1],y])\right\rangle \right) \right| _{t=0}
\\
=-2\left\langle y,[[y,\theta y_1],\theta y]\right\rangle \phi ^{\prime
\prime }\left( -\left\langle y,\theta y\right\rangle \right) 
\end{eqnarray*}
Since $y$ is a ${\frak k}$-conjugate to a root vector, there is a scalar $%
k^{\prime }$ independent of $y$ such that $[[y,\theta y],y]=k^{\prime
}\left\langle y,\theta y\right\rangle y$ . Setting $y=y_1$ we get $%
\left\langle y_1,\theta y_1\right\rangle =-1$, 
\[
[[y_1,-x_1],y_1]=-2y_1
\]
and $k^{\prime }=2.$ Also $-\left\langle y,[[y,\theta y_1],\theta
y]\right\rangle =\left\langle [y,\theta y],[y,\theta y_1]\right\rangle
=\left\langle [[y,\theta y],y],\theta y_1\right\rangle $. Hence 
\begin{equation}
\sum_jc_jl_j\cdot \phi ^{\prime }=4\left\langle y,\theta y\right\rangle
\left\langle \theta y_1,y\right\rangle \phi ^{\prime \prime }\text{.}
\label{=term2}
\end{equation}

\item  Next we note that $\sum_jl_j\cdot c_j$ is independent of the basis $%
l_j$, so we may assume that 
\[
\theta l_j=\pm l_j\text{ and }\left\langle l_j,-\theta l_k\right\rangle
=\delta _{jk}\text{.}
\]
Then $c_j(y)=\left\langle [\theta y,y_1],-\theta l_j\right\rangle $ and 
\begin{eqnarray*}
\tsum\nolimits_jl_j\cdot c_j=\tsum\nolimits_j\left\langle [\theta
[l_j,y],y_1],-\theta l_j\right\rangle  \\
=\tsum\nolimits_j\left\langle y_1,[\theta [l_j,y],\theta l_j]\right\rangle
=-\left\langle y_1,\Omega \theta y\right\rangle .
\end{eqnarray*}
Here $\Omega =\tsum\nolimits_j\limfunc{ad}(\theta l_j)^2=\Omega _{{\frak l}%
}-2\Omega _{{\frak k}\cap {\frak l}}\,$, where the Casimir elements are
obtained by using dual bases with respect to $\left\langle
\,,\,\right\rangle $.

To continue, we need the following lemma:
\end{itemize}

\begin{lemma}
$\Omega $ acts on ${\frak n}$ by the scalar $k^{\prime \prime }=2-2e$.\label
{kpp}
\end{lemma}

\TeXButton{Proof}{\proof} When $e=1$ it's easy to see that the operator $%
\Omega $ acts by $0$. Indeed, \ in this case ${\frak g}$ is a complex
semisimple Lie algebra and for each basis element $l_j\in {\frak k}\cap 
{\frak l}$ there exists a basis element $l_j^{\prime }=\sqrt{-1}l_j\in 
{\frak p}\cap {\frak l}$ . Then $[l_j,[l_j,x]]+[l_j^{\prime },[l_j^{\prime
},x]]=0$ and $k^{\prime \prime }=0$.

When $e=0,$ ${\frak g}$ is split and simply laced, and ${\frak l}$ is the
split real form of a complex reductive algebra ${\frak l}_{{\Bbb C}}$. Take
a root vector $x_\lambda \in {\frak g}_\lambda $, where $\lambda $ is any
positive root in ${\frak n}$. For any positive root $\alpha $ of ${\frak l}_{%
{\Bbb C}}$ we fix $e_\alpha \in {\frak l}_\alpha $ and set $l_\alpha
=e_\alpha +\theta e_\alpha \in {\frak k}\cap {\frak l}$ and $l_\alpha
^{\prime }=e_\alpha -\theta e_\alpha \in {\frak p}\cap {\frak l}$. Then the
collection of all $l_\alpha $, $l_\alpha ^{\prime }$ together with the
orthonormal basis of a Cartan subalgebra ${\frak f}$ of ${\frak l}$ forms a
basis of ${\frak l}$. Observe that 
\[
\lbrack l_\alpha ,[l_\alpha ,x_\lambda ]]+[l_\alpha ^{\prime },[l_\alpha
^{\prime },x_\lambda ]]=[e_\alpha ,[e_\alpha ,x_\lambda ]]+[e_{-\alpha
},[e_{-\alpha },x_\lambda ]]=0\text{,} 
\]
since $x_\lambda \in {\frak g}_\lambda $ and neither $\lambda +2\alpha $ nor 
$\lambda -2\alpha $ is a root of the simply laced algebra ${\frak g}_{{\Bbb C%
}}$.

We choose a basis $\{u_i\}$ of ${\frak f}$, and denote the elements of the
dual (with respect to $\left\langle \;,\;\right\rangle $) basis by $%
\widetilde{u}_i.$ Then 
\[
\Omega x_\lambda =\sum_i[u_i,[\widetilde{u}_i,x_\lambda ]]=\left\langle
\lambda ,\lambda \right\rangle x_\lambda =2x_\lambda \text{.} 
\]

In the remaining two cases ${\frak k}\cap {\frak l}$ acts on ${\frak n}$
irreducibly, therefore $\Omega $ automatically acts by a scalar and it
suffices to compute $\sum_j[l_j,[l_j,x_1]]$. For $e=3$ we have $G=GL_{2n}(%
{\Bbb H}),L=GL_n({\Bbb H})\times GL_n({\Bbb H})$ and ${\frak n}={\Bbb H}%
^{n\times n}$. The computation for this group is similar to the case of $%
G=GL_{2n}({\Bbb R})$. We reduce the calculation to the summation over the
diagonal subalgebra of ${\frak l}$ and obtain

\[
\Omega x_\lambda =\left\langle \lambda ,\lambda \right\rangle x_\lambda
+3\left\langle \sqrt{-1}\lambda ,\sqrt{-1}\lambda \right\rangle x_\lambda
=-4x_\lambda \text{.} 
\]

Finally, for $e=2$ $(G=Sp_{n,n})$, a direct evaluation of $%
\sum_j[l_j,[l_j,x_1]]$ gives $k^{\prime \prime }=-2$.\TeXButton{End Proof}
{\endproof} \smallskip\ \\\ Therefore, we get 
\begin{equation}
\sum_j\phi ^{\prime }l_j\cdot c_j=-2(1-e)\left\langle \theta
y_1,y\right\rangle \phi ^{\prime }\text{.}  \label{=term3}
\end{equation}

\begin{itemize}
\item  Finally, we have 
\begin{equation}
\pi (\theta y_1)f=\frac d{dt}\left. \left( e^{-i\left\langle x+t\theta
y_1,y\right\rangle },\phi \right) \right| _{t=0}=-i\left( e^{-i\left\langle
x,y\right\rangle },\left\langle \theta y_1,y\right\rangle \phi \right) .
\label{=term4}
\end{equation}
\smallskip\ 
\end{itemize}

Putting the formulas (\ref{=term1})--(\ref{=term4}) together, we deduce the
lemma.\TeXButton{End Proof}{\endproof}\smallskip\ 

{\bf Proof of Proposition \ref{=spher}.} Recall that we study $\phi _\tau
(z)=\dfrac{K_\tau (\sqrt{z})}{\left( \sqrt{z}\right) ^\tau }$ , its lift $%
g_\tau $ to the radial function on ${\cal O}_1$, 
\[
g_\tau (y)=\phi _\tau (-\left\langle y,\theta y\right\rangle )=\dfrac{K_\tau
(\left| y\right| )}{\left| y\right| ^\tau } 
\]
and its Fourier transform $\Phi (x)=(e^{-i\left\langle x,y\right\rangle
},g_\tau )$. By Proposition \ref{=infin} it suffices to check that $\pi
(y_1+\theta y_1)\Phi =0$. This identity would follow immediately from Lemma 
\ref{=main-mess}, because $D\phi _\tau =0$ by formula (\ref{=change}) and
then the desired result follows from (\ref{=crown}).

To complete the proof we have to verify the assumptions (\ref{=diff-cond}).
In subsection \ref{=analytic} we proved that $g_\tau \in L^1({\cal O}_1,d\mu
_1)$. It is easy to verify (using the standard facts about the derivatives
of $K_\tau $ from \cite{watson}), that the lifts to ${\cal O}_1$ of the
functions $\phi _\tau ^{\prime }(z)$ and $\phi _\tau ^{\prime \prime }(z)$
(we denote them by $g_\tau ^{\prime }(y)$ and $g_\tau ^{\prime \prime }(y)$)
both belong to $L^1({\cal O}_1,d\mu _1)$. Observe also that $\phi _\tau (z)$%
, $\phi _\tau ^{\prime }(z)$, $\phi _\tau ^{\prime \prime }(z)$ are all
monotone on $(0,\infty )$.

Moreover, since all these functions tend to zero exponentially as $\left|
y\right| \rightarrow \infty $, the functions $A(y)g_\tau (y)$, $A(y)g_\tau
^{\prime }(y)$, $A(y)g_\tau ^{\prime \prime }(y)$ all belong to $L^1({\cal O}%
_1,d\mu _1)$, for any $A(y)$ bounded in the neighbourhood of $y=0$ and
growing (at most) polynomially with respect to $\left| y\right| $ as $\left|
y\right| \rightarrow \infty $.

Fix $h\in {\frak l}$ , $x\in {\frak n}$ and choose $c>0$ sufficiently small,
such that for all $y\in {\cal O}_1$ and $\left| t\right| <c$%
\[
\left| \left\langle e^{th}\cdot y,\theta (e^{th}\cdot y)\right\rangle
\right| \geq \frac{\left| \left\langle y,\theta y\right\rangle \right| }2%
\text{.} 
\]
We can then estimate the derivative: 
\[
\left| \frac d{dt}\left( e^{-i\left\langle x,e^{th}\cdot y\right\rangle
}\phi _\tau \left( -\left\langle e^{th}y,\theta e^{th}y\right\rangle \right)
\right) \right| \leq \left| A_1(y)\phi _\tau \left( \left| y\right|
^2/2\right) \right| +\left| A_2(y)\phi _\tau ^{\prime }\left( \left|
y\right| ^2/2\right) \right| \text{,} 
\]
for all $y\in {\cal O}_1$ and $\left| t\right| <c$, where $A_1(y),A_2(y)$
are some functions of polynomial growth. From the discussion above, the
right-hand side of this inequality is an $L^1$-function on ${\cal O}_1$,
hence $e^{-i\left\langle x,y\right\rangle }g_\tau \in {\cal I}$.

Proceeding in the same manner, we deduce that $e^{-i\left\langle
x,y\right\rangle }g_\tau ^{\prime }\in {\cal I}$. \TeXButton{End Proof}
{\endproof} \smallskip\ 

\subsection{Proof of Theorem 0.1}

Denote by ${\bf J}$ the space of the induced representation $\pi _1=\limfunc{%
Ind}_{\overline{P}}^G(e^{-d\nu })$. By the Gelfand-Naimark decomposition and
the $\exp $ map, ${\bf J}$ can be viewed as a subspace of $C^\infty ({\frak n%
})$. Then for $l\in L$ and $\eta \in {\bf J}$ we have 
\[
\pi _1(l)\eta (x)=e^{-d\nu }(l)\eta (l^{-1}\cdot x)\text{.} 
\]
It was proved in \cite{sahi-dp} that the $({\frak g},K)$-module ${\bf J}$
has a unitarizable spherical $({\frak g},K)$-submodule $V$, which we also
regard as a subspace of $C^\infty ({\frak n})$.

{\bf Remark. }It is possible to give a direct description of the elements of
the ``abstract'' Hilbert space ${\cal H}$, where ${\cal H}$ is the Hilbert
space closure of $V$ with respect to the $({\frak g},K)$-invariant norm on $%
V $. For that purpose we use the ``compact'' realization of $\pi _1$ on $%
C^\infty (K/M)$ from \cite{sahi-dp}. It was shown that $\pi _1$ is a
representation of ladder type, with all its $K$-types $\{\alpha _m\,|\;m\in 
{\Bbb N}\}$ lying on a single line, $\alpha _1$ being a one-dimensional $K$%
-type. The restriction $\left\langle \;,\;\right\rangle _m$ of a $\pi _1$%
-invariant Hermitian form to any $K$-type $\alpha _m$ is a multiple of the $%
L^2(K)$-inner product on $V$, and from the explicit formulas in \cite
{sahi-dp} it follows that 
\[
q_m\stackrel{\text{def}}{=}\frac{\left\langle \;,\;\right\rangle _m}{%
\left\langle \;,\;\right\rangle _1}=O(m^C) 
\]
for some constant $C>1$, which can be expressed in terms of parameters $d,$ $%
e$ and $n$. Thus we can identify ${\cal H}$ with the Hilbert space $%
L^2\left( {\Bbb N},\{q_m\}\right) $, where the constant $q_m$ gives the
weight of the point $m\in {\Bbb N}$. That is, any element of ${\cal H}$ can
be viewed as an $M$-equivariant function on $K$, such that its sequence of
Fourier coefficients belongs to $L^2\left( {\Bbb N},\{q_m\}\right) $. In
particular $L^2\left( {\Bbb N},\{q_m\}\right) \subset l^2({\Bbb N)}$, and
the elements of ${\cal H}$ all lie in $L^2(K)$.\smallskip\ 

We write ${\bf H}$ for the space of those tempered distributions on ${\frak n%
}$ which are Fourier transforms of $\psi d\mu _1$ for some $\psi \in L^2(%
{\cal O}_1,d\mu _1)$. If $\eta $ is the Fourier transform of a distribution
of the form $\psi d\mu _1$, i.e., 
\[
\eta (x)=\int_{{\cal O}_1}e^{-i\left\langle x,y\right\rangle }\psi (y)d\mu
_1(y)=\left( e^{-i\left\langle x,y\right\rangle },\psi (y)\right) ,
\]
then 
\begin{eqnarray*}
\pi _1(l)\eta (x) &=&e^{-d\nu }(l)\eta (l^{-1}\cdot x)=\int_{{\cal O}%
_1}e^{-i\langle l^{-1}\cdot x,y\rangle }\psi (y)e^{-d\nu }(l)d\mu _1(y) \\
&=&\int_{{\cal O}_1}e^{-i\langle l^{-1}x,l^{-1}y\rangle }\psi (l^{-1}\cdot
y)e^{-dv}(l)d\mu _1(l^{-1}\cdot y) \\
&=&(e^{-i\left\langle x,y\right\rangle },e^{d\nu }(l)\psi (l^{-1}\cdot y)).
\end{eqnarray*}
It follows from the calculation above that $P$ acts {\em unitarily }on{\em \ 
}${\bf H}$ (it is convenient to identify this action with its realization on 
$L^2({\cal O}_1,d\mu _1)$ via the Fourier transform).We denote this unitary
representation of $P$ by $\pi ^{\prime }$. Observe that $(\pi ^{\prime },%
{\bf H})$ is an irreducible representation of $P$.

According to Proposition \ref{=sq-int}, $\Phi (x)=(e^{-i<x,y>},\left|
y\right| ^{-\tau }K_\tau (\left| y\right| ))$ belongs to ${\bf H}$.

\begin{theorem}
$V$ is a dense subspace of ${\bf H}$, and the restriction of the norm is $(%
{\frak g},K)$-invariant.
\end{theorem}

\TeXButton{Proof}{\proof} Let $C^\infty (K)_V$ be the subspace of $C^\infty
(K)$, consisting of those smooth functions on $K$, whose $K$-isotypic
components belong to $V$. Since $V$ is a submodule of ${\bf J}$, $C^\infty
(K)_V$ is obviously $G$-invariant.

Denote by ${\cal C}(G)$ the convolution algebra of smooth $L^1$ functions on 
$G=PK$, and consider 
\[
{\bf W}=\pi _1({\cal C}(G))\Phi \subset C^\infty (K)_V. 
\]
So all elements of ${\bf W}$ are continuous functions on $K$, hence
continuous on $G$, and therefore are determined by their restrictions to $N$%
. Moreover, ${\bf W}=\pi _1({\cal C}(PK))\Phi $ and $K$ fixes $\Phi $,
therefore 
\[
{\bf W}=\pi _1({\cal C}(P))\Phi =\pi ^{\prime }({\cal C}(P))\Phi \text{.} 
\]
This shows that ${\bf W}$ is a $\pi ^{\prime }(P)$-invariant subspace of $%
{\bf H}$, and from the irreducibility of $\pi ^{\prime }$ we conclude that $%
{\bf W}$ is dense in ${\bf H}$.

We can now put two $\pi _1(P)$-invariant norms on ${\bf W}$ -- one from $%
{\bf H}$ and another from $V$, as follows. If $f=\sum c_mv_m$, with $v_m$ in
the $K$-isotypic component with highest weight $\alpha _m$ (occurring in $V$%
) and $\left\| v_m\right\| _{L^2(K)}=1$, then 
\begin{equation}
\left\| f\right\| _V^2=\sum \left| c_m\right| ^2q_m\text{.}
\label{=new-norm}
\end{equation}
Since $f$ \thinspace is smooth, it follows that $\left| c_m\right| $ decays
rapidly, so the series in (\ref{=new-norm}) converges, thus giving a $\pi
_1(P)$-invariant norm on ${\bf W}$.

Then it follows from \cite{poguntke} (cf. \cite[p.417]{sahi-expl}), that we
can find a (dense) ${\cal C}(P)$-invariant subspace ${\bf W}^{\prime
}\subset {\bf W}$, such that these two forms are proportional on ${\bf W}%
^{\prime }$. Considering the closure of ${\bf W}^{\prime }$ we obtain an
isometric $P$-invariant imbedding of ${\bf H}$ into ${\cal H}$.

Then ${\bf W}$ is:

(1) \ a $G$-invariant subspace of the irreducible module ${\cal H}$, hence
dense in ${\cal H}$;

(2) \ a dense subspace of the Hilbert space ${\bf H}$.

It follows that ${\bf H}={\cal H}$. \TeXButton{End Proof}{\endproof}%
\smallskip\ 

This concludes the proof of Theorem \ref{=Main-theorem}.

\section{Tensor powers of $\pi _1$}

\subsection{Restrictions to $P$}

In the previous section we constructed a unitary representation $\pi _1$ of $%
G$ acting on the Hilbert space $L^2({\cal O}_1,d\mu _1)$, where ${\cal O}_1$
is the minimal $L$-orbit in a non-Euclidean Jordan algebra $N$. Define the $%
k $-th tensor power representation 
\[
\Pi _k=\pi _1^{\otimes k}(2\leq k<n). 
\]
As we shall show, the techniques developed in \cite{tens} allow us to
establish a duality between the spectrum of this tensor power and the
spectrum of a certain homogeneous space. We omit the proofs of the several
propositions below, because the proofs of the corresponding statements from 
\cite{tens} can be used without any substantial modification.

Observe that the orbit ${\cal O}_k$ is dense in $\stackunder{k\text{ times}}{%
\underbrace{{\cal O}_1+{\cal O}_1+\ldots +{\cal O}_1}}$. The representation $%
\Pi _k$ acts on $\left[ L^2({\cal O}_1,d\mu _1)\right] ^{\otimes k}\simeq
L^2({\cal O}_k^{\prime },d\mu ^{\prime })$, where ${\cal O}_k^{\prime }=%
{\cal O}_1^{\times k}$ and $d\mu ^{\prime }$ is the product measure on${\cal %
\ O}_k^{\prime }$ . We fix a generic representative $\xi ^{\prime }=(\xi
_1,\xi _2,\ldots ,\xi _k)\in {\cal O}_k^{\prime }$ , such that 
\[
\xi =\xi _1+\xi _2+\ldots +\xi _k\in {\cal O}_k. 
\]

Denote by $S_k^{\prime }$ and $S_k$ the isotropy subgroups of $\xi ^{\prime
} $ and $\xi $, respectively, with respect to the action of $L$ on ${\cal O}%
_k^{\prime }$ and ${\cal O}_k$. Observe that the Lie algebras ${\frak s}%
_k^{\prime }$ and ${\frak s}_k$ of $S_k^{\prime }$ and $S_k$, respectively,
can be written as 
\begin{eqnarray*}
{\frak s}_k^{\prime } &=&({\frak h}_k+{\frak l}_k)+{\frak u}_k \\
{\frak s}_k &=&({\frak g}_k+{\frak l}_k)+{\frak u}_k\text{.}
\end{eqnarray*}
Here ${\frak l}_k,{\frak g}_k$ and ${\frak h}_k$ are reductive, ${\frak h}%
_k\subset {\frak g}_k$ and ${\frak u}_k$ is a nilpotent radical common for
both ${\frak s}_k^{\prime }$ and ${\frak s}_k$. Let $G_k$ and $H_k$ be the
corresponding Lie groups.\smallskip\ 

{\bf Example. }Take $G=O_{2n,2n}$ and $k<n$. Then $\xi
_i=E_{2i-1,2i}-E_{2i,2i-1}$ ($1\leq i\leq k$), $\xi =\sum_{i=1}^k\xi _i$ and 
\[
{\frak s}_k=\left( {\frak sp}_{2k}({\Bbb R})+{\frak gl}_{2(n-k)}({\Bbb R)}%
\right) +{\Bbb R}^{2k,2(n-k)}\text{.} 
\]
Then $G_k=Sp_{2k}({\Bbb R)}$ and it's easy to check that $H_k=SL_2({\Bbb R}%
)^k$.\smallskip\ 

The following Lemma can be verified by direct calculation (cf. 
\cite[Lemma 2.1]{tens}).

\begin{lemma}
Let $\chi _\xi $ be the character of $N$ corresponding to $\xi \in N^{*}$.
Then 
\[
\Pi _k|_P=\limfunc{Ind}\nolimits_{S_k^{\prime }N}^P(1\otimes \chi _\xi )=%
\limfunc{Ind}\nolimits_{S_kN}^P\left( (\limfunc{Ind}\nolimits_{S_k^{\prime
}}^{S_k}1)\otimes \chi _\xi \right) \;\;\;\text{{\em (}}L^2\text{{\em %
-induction).}\TeXButton{End Proof}{\endproof}} 
\]
\end{lemma}

Let $\gamma ^{\prime }=\limfunc{Ind}_{H_k}^{G_k}1$ be the quasiregular
representation of $G_k$ on $L^2(G_k/H_k)$; then it can be decomposed using
the Plancerel measure $d\mu $ for the reductive homogeneous space $%
X_k=G_k/H_k$ and the corresponding multiplicity function $m:\widehat{G}%
_k\rightarrow \{0,1,2,\ldots \}$, i.e., 
\[
\gamma ^{\prime }\simeq \int_{\widehat{G}_k}^{\oplus }m(\kappa )\kappa
\,d\mu (\kappa )\text{.} 
\]
Each irreducible representation $\kappa $ of $G_k$ can be extended to an
irreducible representation $\kappa ^{\vee }$ of $S_k$, and the decomposition
of the Lemma above can be rewritten as 
\begin{equation}
\Pi _k|_P=\int_{\widehat{G}_k}^{\oplus }m(\kappa )\Theta (\kappa )\,d\mu
(\kappa )\text{,}  \label{=Pi-to-P}
\end{equation}
where $\Theta (\kappa )=\limfunc{Ind}_{S_kN}^P(\kappa ^{\vee }\otimes \chi
_\xi )$.

Moreover, by Mackey theory all representations $\Theta (\kappa )$ are
unitary irreducible representations of $P$, and $\Theta (\kappa ^{\prime
})\simeq \Theta (\kappa ^{^{\prime \prime }})$ if and only if $\kappa
^{\prime }\simeq \kappa ^{\prime \prime }$.

\subsection{Low-rank theory}

In \cite{tens} we extended the theory of low-rank representations (\cite
{li-singular}) to the conformal groups of euclidean Jordan algebras.
Inspection of the argument in \cite{tens} shows that the analogous theory
can be developed in exactly the same manner for the conformal groups of {\em %
non-euclidean} Jordan algebras.

For any unitary representation $\eta $ of $G$, we decompose its restriction $%
\eta |_N$ into a direct integral of unitary characters, where the
decomposition is determined by a projection-valued measure on $\widehat{N}%
=N^{*}$. If this measure is supported on a single {\em non-open} $L$%
-orbit\thinspace {\em \ }${\cal O}_m$,\thinspace $1\leq m<n$ we call $\eta $
a {\em low-rank representation}, and write $\limfunc{rank}\eta =m$.
Proceeding by induction on $m$, as in \cite{li-singular}, \cite[Sect 3]{tens}%
, we can prove the following

\begin{theorem}
Let $\eta $ be a low-rank representation of $G.$ Write ${\cal A}(\eta ,P)$
for the von Neumann algebra generated by $\{\eta (x)|\,x\in P\}$ and ${\cal A%
}(\eta ,G)$ for the von Neumann algebra generated by $\{\eta (x)|\,x\in G\}$%
. Then ${\cal A}(\eta ,G)={\cal A}(\eta ,P)$. \TeXButton{End Proof}
{\endproof}\smallskip\ 
\end{theorem}

{\bf Proof of Theorem \ref{=Main-tensoring}. }Now consider the restriction
of $\Pi _k$ to $N$. Its restriction to $P$ is given by the direct integral
decomposition (\ref{=Pi-to-P}), and we can further restrict it to $N$. The
rank of the induced representation $\Theta (\kappa )=\limfunc{Ind}%
_{S_kN}^P(\kappa ^{\vee }\otimes \chi _\xi )$ is $k$ (the $N$--spectrum is
supported on the $L-$orbit of $\xi $, i.e. ${\cal O}_k$). Therefore $\Pi _k$
can be decomposed over the irreducible representations of $G$ of rank $k$.

It follows from the theorem above that any two non-isomorphic
representations from the spectrum of $\Pi _k$ restrict to non-isomorphic
irreducible representations of $P$. Hence the representation $\Pi _k$ can be
decomposed as 
\begin{equation}
\Pi _k=\int_{\widehat{G}_k}^{\oplus }m(\kappa )\theta (\kappa )\,d\mu
(\kappa ),  \label{=Pi-decomp}
\end{equation}
where for almost every $\kappa $ the unitary irreducible representation $%
\theta (\kappa )$ is obtained as the {\em unique }irreducible representation
of $G$ determined by the condition $\theta (\kappa )|_P=\Theta (\kappa )$.

Therefore, the map $\kappa \rightarrow \theta (\kappa )$ gives a
(measurable) bijection between the spectrum of $\Pi _k=\pi ^{\otimes k}$ and
the unitary representations of $G_k$ occurring in the quasiregular
representation on $L^2(G_k/H_k).$\TeXButton{End Proof}{\endproof}\smallskip\ 

{\bf Example. }Take $G=E_{7(7)}$. It is the conformal group of the split
exceptional real Jordan algebra $N$ of dimension 27. Consider the tensor
square of the minimal representation $\pi _1$ of $G$ ($k=2$). Then $L={\Bbb R%
}^{*}\times E_{6(6)}$, $S_2^{\prime }$ is the stabilizer of $\,y_1$ {\bf and 
}$y_2$ and $S_2$ is the stabilizer of $y_1+y_2\in {\cal O}_2$. One can see
that in this case ${\frak g}_2=\limfunc{Stab}_{{\frak s0}(5,5)}(y_1+y_2)=%
{\frak so}(4,5)$ and ${\frak h}_2=\limfunc{Stab}_{{\frak s0}(5,5)}(y_1)\cap 
\limfunc{Stab}_{{\frak s0}(5,5)}(y_2)={\frak so}(4,4)$ (cf. \cite[16.7]
{adams}). Hence the decomposition (\ref{=Pi-decomp}) establishes a duality
between the representations of $E_{7(7)}$ occurring in $\Pi _2=\pi _1\otimes
\pi _1$ and the unitary representations of $Spin(4,5)$ occurring in $%
L^2\left( Spin(4,5)/Spin(4,4)\right) $ . The homogeneous space $%
Spin(4,5)/Spin(4,4)$ is a (pseudo-riemannian) symmetric space of rank 1, and
it is known to be multiplicity free. Therefore, $\pi _1\otimes \pi _1$ has
simple spectrum.

Similarly, for $G=E_7({\Bbb C})$ we obtain a duality between $E_7({\Bbb C})$
and the symmetric space $SO_9({\Bbb C})/SO_8({\Bbb C})$.\newpage\ 

\ \appendix 

\section{Groups associated to non-Euclidean Jordan algebras\protect\medskip\ 
}

\begin{tabular}{|c|l|c|c|c|}
\hline
$G$ & $K/M$ & $d$ & $e$ & $G_k/H_k$ for $2\leq k<n$ \\ \hline\hline
\multicolumn{1}{|l|}{$GL_{2n}({\Bbb R})$} & $O_{2n}/(O_n\times O_n)$ & 
\multicolumn{1}{|l|}{$1$} & \multicolumn{1}{|l|}{$0$} & \multicolumn{1}{|l|}{%
$GL_k({\Bbb R})/[GL_1({\Bbb R})]^k$} \\ \hline
\multicolumn{1}{|l|}{$O_{2n,2n}$} & $(O_{2n}\times O_{2n})/O_{2n}$ & 
\multicolumn{1}{|l|}{$2$} & \multicolumn{1}{|l|}{$0$} & \multicolumn{1}{|l|}{%
$Sp_{2k}({\Bbb R})/[SL_2({\Bbb R})]^k$} \\ \hline
\multicolumn{1}{|l|}{$E_{7(7)}$} & $SU_8/Sp_4$ & \multicolumn{1}{|l|}{$4$} & 
\multicolumn{1}{|l|}{$0$} & \multicolumn{1}{|l|}{$Spin(4,5)/Spin(4,4)$ \ }
\\ \hline
\multicolumn{1}{|l|}{$O_{p+2,p+2}$} & $[O_{p+2}]^2/{\Bbb [}O_1\times
O_{p+1}^2]$ & \multicolumn{1}{|l|}{$p$} & \multicolumn{1}{|l|}{$0$} & 
\multicolumn{1}{|l|}{} \\ \hline\hline
\multicolumn{1}{|l|}{$Sp_n({\Bbb C})$} & $Sp_n/U_n$ & \multicolumn{1}{|l|}{$%
1 $} & \multicolumn{1}{|l|}{$1$} & \multicolumn{1}{|l|}{$O_k({\Bbb C})/[O_1(%
{\Bbb C})]^k$} \\ \hline
\multicolumn{1}{|l|}{$GL_{2n}({\Bbb C})$} & $U_{2n}/(U_n\times U_n)$ & 
\multicolumn{1}{|l|}{$2$} & \multicolumn{1}{|l|}{$1$} & \multicolumn{1}{|l|}{%
$GL_k({\Bbb C})/[GL_1({\Bbb C})]^k$} \\ \hline
\multicolumn{1}{|l|}{$O_{4n}({\Bbb C})$} & $O_{4n}/U_{2n}$ & 
\multicolumn{1}{|l|}{$4$} & \multicolumn{1}{|l|}{$1$} & \multicolumn{1}{|l|}{%
$Sp_{2k}({\Bbb C})/[SL_2({\Bbb C})]^k$} \\ \hline
\multicolumn{1}{|l|}{$E_7({\Bbb C})$} & $E_7/(E_6\times U_1)$ & 
\multicolumn{1}{|l|}{$8$} & \multicolumn{1}{|l|}{$1$} & \multicolumn{1}{|l|}{%
$SO_9({\Bbb C})/SO_8({\Bbb C})$ \ } \\ \hline
\multicolumn{1}{|l|}{$O_{p+4}({\Bbb C})$} & $O_{p+4}/(O_{p+2}\times U_1)$ & 
\multicolumn{1}{|l|}{$p$} & \multicolumn{1}{|l|}{$1$} & \multicolumn{1}{|l|}{
} \\ \hline\hline
\multicolumn{1}{|l|}{$Sp_{n,n}$} & $(Sp_n\times Sp_n)/Sp_n$ & 
\multicolumn{1}{|l|}{$2$} & \multicolumn{1}{|l|}{$2$} & \multicolumn{1}{|l|}{%
$O_k^{*}/[O_1^{*}]^k$} \\ \hline
\multicolumn{1}{|l|}{$GL_{2n}({\Bbb H)}$} & $Sp_{2n}/(Sp_n\times Sp_n)$ & 
\multicolumn{1}{|l|}{$4$} & \multicolumn{1}{|l|}{$3$} & \multicolumn{1}{|l|}{%
$GL_k({\Bbb H})/[GL_1({\Bbb H})]^k$} \\ \hline
\end{tabular}
\medskip\

\end{document}